\title{ ~~\\ Addendum to the paper: ``Artin Prime Producing Quadratics''  [Abh. Math. Sem. Univ. Hamburg 77 (2007), 109--127; 
MR2379332 (2008m:11194)] by P. Moree}
\author{Yves Gallot and Pieter Moree}
\documentclass[12pt]{article}
\usepackage{amssymb, latexsym, amsfonts, makeidx}
\usepackage{amsmath}
\def\@ptsize{2}
\begin{document}
\date{}
\maketitle 
{\def\thefootnote{}
\begin{abstract}
\noindent A record mentioned in the paper was recently improved on by Akbary and Scholten. However, the
record mentioned was not the then record. The then record, due to Gallot (2004), actually slightly improves on
that obtained recently by Akbary and Scholten.
\end{abstract}
Given an integer $g$ and a polynomial $f(X)\in \Bbb Z[X]$, let $p_1(g,f),p_2(g,f),\ldots$ be
the sequence of primes that is obtained on going through the sequence $f(0),f(1),\ldots$ and writing
down the primes not dividing $g$ as they appear (called Artin primes). We let $r$ be the largest integer $r$ (if this
exists) such that $g$ is a primitive root mod $p$ for all primes $p_j(g,f)$ with $1\le j\le r$. We let
$c_g(f)$ be the number of distinct primes amongst $p_j(g,f)$ with $1\le j\le r$.

In \cite{Moree} the problem was addressed of finding an integer $g$ and a quadratic polynomial $f$ such
that $c_g(f)$ is as large as possible and it was stated that $$c_g(f)=31082$$ was the current record (obtained
by Yves Gallot).
On preparing the paper for publication (fall 2006) the author failed to recall an e-mail by Gallot from June
2004. That e-mail actually stated what in  2006 still would be the true current record (due to Gallot), namely $$c_g(f)=38639.$$

It is obtained on taking $f(X)=32 X^2 + 39721664 X + 182215381147285848449$
and  $g = 593856338459898$. Perhaps a more elegant reformulation is:
for those $38639$ integers $n$ in $[620651, 1749283]$ for which $$h(n):=32 n^2 + 182215368820640606817$$ is prime, 
the number $593856338459898$ is a primitive root modulo $h(n)$.

In a recent paper by Akbary and Scholten \cite{AS} the authors find a $g$ and a quadratic $f$ such that
$c_g(f)=37951$. This improves on the record indicated in \cite{Moree}, but falls slightly below the `hidden record' 
indicated above.

Akbary and Scholten go beyond Moree in that they in addition consider the case where $f$ linear and $f$ cubic 
and obtain here record values for consecutive Artin primes for certain integers $g$ of $6355$, respectively $10011$.

Finally, let us mention some highly interesting work by Pollack \cite{P}. 
He merges the method of proof of Hooley \cite{Hooley} of Artin's conjecture (under GRH) with the method
of Maynard-Tao \cite{mei, polymath} in order to produce bounded gaps between primes: On GRH for every nonsquare $g\ne -1$ and
every $m$, there are infinitely many runs of $m$ consecutive primes all possessing $g$ as a primitive root
and lying in an interval of order $O_m(1)$. For related work see Baker and Pollack \cite{BP}.

\end{document}